\documentclass[12pt]{amsart}
\usepackage{amsfonts}
\usepackage{amssymb}
\usepackage{latexsym}

\newtheorem{theorem}{Theorem}[section]
\newtheorem{corollary}[theorem]{Corollary}
\newtheorem{lemma}[theorem]{Lemma}
\newtheorem{proposition}[theorem]{Proposition}

\newtheorem{problem}[theorem]{Problem}

\theoremstyle{definition}
\newtheorem{definition}[theorem]{Definition}
\newtheorem{remark}[theorem]{Remark}
\newtheorem{example}[theorem]{Example}

\newtheorem{algorithm}[theorem]{Algorithm}

\begin{document}

\title{Gomory integer programs}
\author{Serkan Ho{\c s}ten}
\address{Department of Mathematics, San Francisco State University, 
San Francisco, CA 94132}
\email{serkan@math.sfsu.edu}
\author{Rekha R. Thomas}
\address{Department of Mathematics, University of Washington,
Box 354350, Seattle, WA 98195-4350}
\email{thomas@math.washington.edu}

\date{\today}
\begin{abstract}
  The set of all group relaxations of an integer program contains
  certain special members called Gomory relaxations. A family of
  integer programs with a fixed coefficient matrix and cost vector but
  varying right hand sides is a Gomory family if every program in the
  family can be solved by one of its Gomory relaxations. In this
  paper, we characterize Gomory families. Every TDI system gives a
  Gomory family, and we construct Gomory families from matrices whose
  columns form a Hilbert basis for the cone they generate. The
  existence of Gomory families is related to the Hilbert covering
  problems that arose from the conjectures of Seb{\"o}. Connections to
  commutative algebra are outlined at the end.
\end{abstract}

\maketitle

\section{Introduction}
Given an integer $d \times n$ matrix $A$ and a cost vector $c \in
\mathbb Z^n$, we consider the family $IP_{A,c}$ of all feasible
integer programs of the form
$$IP_{A,c}(b) := min \,\,\{ c \cdot x \,\, : \,\, Ax = b, \,\, x \in
\mathbb N^n \}$$ as the right hand side vector $b$ varies. 
The
matrix $A$ is assumed to have rank $d$ and $cone(A)$, the cone 
generated by the columns of $A$, is assumed to be pointed.
We also assume that $\{x \in \mathbb R^n_{\geq 0} \,
\, : \,\, Ax = 0 \} = \{0\}$ which guarantees that all programs in
$IP_{A,c}$ are bounded.

In \cite{Gom65}, Gomory defined the {\em group
relaxation} of $IP_{A,c}(b)$: 
$$
min \, \, \{ \tilde{c}_{{\bar \sigma}} \cdot x_{\bar \sigma} \,\, :
A_{\sigma}x_{\sigma} + A_{\bar \sigma}x_{\bar \sigma} = b, \,\,
x_{\bar \sigma} \geq 0, (x_{\sigma}, x_{\bar \sigma}) \in \mathbb Z^n
\}$$
where $A_{\sigma}$ is the optimal basis of the linear relaxation
of $IP_{A,c}(b)$ and non-negativity restrictions on the optimal basic
variables $x_{\sigma}$ have been dropped.  The cost vector
$\tilde{c}_{{\bar \sigma}}$ is the restriction of
$(c-c_{\sigma}A_{\sigma}^{-1}A)$ to the components indexed by the
complement of $\sigma$. The {\em extended group relaxations} of
$IP_{A,c}(b)$ introduced by Wolsey are the $2^{|\sigma|}$ relaxations
obtained by dropping non-negativity restrictions on each subset of the
variables in $x_{\sigma}$ \cite{Wol}. The set of $d$-dimensional
simplicial cones $cone(A_{\sigma})$, as $A_{\sigma}$ varies over the
optimal bases of the LP-relaxations of all programs in $IP_{A,c}$,
triangulates $cone(A)$. This is called the {\em regular triangulation}
of $cone(A)$ with respect to $c$ and is denoted as $\Delta_c$. The
collection of all sets $\sigma$ indexing the full dimensional cones of
$\Delta_c$ along with all their subsets form a simplicial complex on
$\{1, \ldots, n\}$ which we also call $\Delta_c$. In this paper we
consider the set of all group relaxations of $IP_{A,c}(b)$ obtained by
dropping non-negativity restrictions on the variables indexed by each
face of $\Delta_c$ (Definition~\ref{grouprels}). This is a larger set
of group relaxations for $IP_{A,c}(b)$ than the set of extended group
relaxations of Wolsey. We show that these are precisely all the
bounded group relaxations of $IP_{A,c}(b)$ (Theorem~\ref{bounded}).
Among these group relaxations, the easiest to solve are those indexed
by the maximal faces of $\Delta_c$. We call these the {\em Gomory
  relaxations} of $IP_{A,c}(b)$. The family of integer programs
$IP_{A,c}$ is called a {\em Gomory family} if all its members can be
solved by one of their Gomory relaxations.

Theorem~\ref{equivconds} characterizes Gomory families. This theorem
is a consequence of re-casting algebraic results on {\em toric initial
  ideals} in terms of group relaxations of integer programs
\cite{HT2}, \cite{HT1}. These algebraic results come from {\em
  Gr\"obner bases} methods in integer programming \cite{Stu}. No
familiarity with these techniques is assumed in this paper. In
Section~3 we relate Gomory families to {\em total dual integrality}
(TDI-ness).  Theorem~\ref{TDI-ness} shows that $yA \leq c$ is a TDI
system if and only if the regular triangulation $\Delta_c$ is {\em
  unimodular}. This leads to Corollary~\ref{TDI-Gomory} that if $yA
\leq c$ is TDI then $IP_{A,c}$ is a Gomory family.

In Sections 4 and 5 we exhibit general classes of Gomory families. A
matrix $A$ is said to be {\em normal} if its columns form a Hilbert
basis for $cone(A)$. In Section~4 we introduce $\Delta$-normal
matrices which form a proper subset of normal matrices.
Theorem~\ref{specialinitial} shows that every $\Delta$-normal matrix
$A$ gives rise to a Gomory family $IP_{A,c}$. While we do not know if
every normal matrix gives rise to a Gomory family, we show that for
small values of $d$ every regular triangulation of $cone(A)$ supports
Gomory families (Theorem~\ref{smalld}). Gomory families induce special
covers of $\mathbb N A$, the semigroup generated by the columns of
$A$, which relates their existence to the Hilbert cover questions
found in \cite{BG}, \cite{FZ} and \cite{Seb}.

Throughout this paper we consider triangulations of $cone(A)$. 
We wish to point out that
a one-dimensional face of any such triangulation must be generated
by a column of $A$. All computations in this paper 
rely on the connections of this material to commutative algebra as
described in \cite{HT2}, \cite{HT1}, \cite{Stu} and \cite{STV}. 
The relevant connections and codes 
are described briefly in Section 6.

\section{Gomory families}
In this paper, we fix a matrix $A \in \mathbb Z^{d
\times n}$ of rank $d$ and a cost vector $c \in \mathbb Z^n$ and consider the
family $IP_{A,c}$ of all integer programs 
$$IP_{A,c}(b) := min \,\,\{ c \cdot x \,\, : \,\, Ax = b, \,\, x \in
\mathbb N^n \}$$ as $b$ varies in the semigroup $\mathbb N A := \{ Au
\, : \, u \in \mathbb N^n \} \subseteq \mathbb Z^d$. This semigroup is
contained in the intersection of $cone(A) := \{ Ax \, : \, x \in
\mathbb R^n_{\geq 0}\}$, and $\mathbb Z A := \{Az \, : \, z \in
\mathbb{Z}^n\}$, the lattice generated by the columns of $A$. We may
assume without loss of generality that $\mathbb ZA = \mathbb Z^d$.

The feasible {\em linear programs} from $A$ and $c$ are of the form 
$$LP_{A,c}(b) := min \,\,\{\,\, c \cdot x \,\, : \,\, Ax = b, \,\, x
\geq 0 \}$$ where $b \in cone(A)$. We denote this family as 
$LP_{A,c}$. For $\sigma \subseteq \{1, \ldots, n\}$, let
$A_{\sigma}$ be the submatrix of $A$ whose set of column indices is
$\sigma$.

\begin{definition} \label{regtriang}
For $\sigma \subseteq \{1, \ldots, n \}$,
$cone(A_{\sigma})$ is a face of the {\em regular subdivision} $\Delta_c$ 
of $cone(A)$ if and only if there exists a vector $y \in \mathbb R^d$
such that $y \cdot a_j = c_j$ for all $j \in \sigma$ and $y \cdot a_j
< c_j$ for all $j \not \in \sigma$.
\end{definition}

\begin{remark}
The regular subdivision $\Delta_c$ is gotten by taking the 
cone in $\mathbb R^{d+1}$ generated by the lifted vectors $(a_i,c_i)
\in \mathbb R^{d+1}$ where $a_i$ is the $i$-th column
of $A$ and $c_i$ is the $i$-th component of $c$, and then projecting the
{\em lower} facets of this lifted cone back onto $cone(A)$. (See
\cite{BFS}.)
\end{remark}

We assume that $c$ is {\em generic}, which means that 
$\Delta_c$ is a {\em triangulation} of $cone(A)$. 
All cost vectors except those lying on a finite set of hyperplanes of
$\mathbb R^n$ are generic \cite{BFS}. Using $\sigma$ to label
$cone(A_{\sigma})$, the triangulation $\Delta_c$ can be denoted as a
set of subsets of $\{1,\ldots, n \}$. This set of sets is closed under
inclusion since $\Delta_c$ is a simplicial complex, and hence it is
specified completely by its maximal elements. 
For a vector $x \in \mathbb R^n$, let $supp(x) := 
\{i : x_i \neq 0 \}$ denote the {\em support} of $x$. The
significance of regular triangulations for linear programming is
summarized in the following proposition.

\begin{proposition} \cite[Lemma 1.4]{ST} \label{optbases}
An optimal solution of $LP_{A,c}(b)$ is any feasible solution 
$x^{\ast}$ such that $supp(x^{\ast}) = \tau$ where $\tau$ is the 
smallest face of the regular triangulation $\Delta_c$ such that $b \in
cone(A_{\tau})$. \qed
\end{proposition}

Proposition~\ref{optbases} implies that $\sigma \subseteq \{1, \ldots,
n\}$ is a maximal face of $\Delta_c$ if and only if $A_{\sigma}$ is an
optimal basis for all $LP_{A,c}(b)$ with $b$ in $cone(A_{\sigma})$.
Given a polyhedron $P \subset \mathbb R^n$ and a face $F$ of $P$, the
{\em normal cone} of $F$ at $P$ is the cone $N_P(F) := \{ \omega \in
\mathbb R^n \, : \, \omega \cdot x' \geq \omega \cdot x, \,\,\forall
x' \in F \text{ and } x \in P\}$. The set of all normal cones of $P$
form the {\it normal fan} of $P$ in $\mathbb R^n$.

\begin{proposition} \label{normalfan}
The regular triangulation $\Delta_c$ of $cone(A)$ is the normal fan of
the polyhedron $P_c := \{ y \in \mathbb R^d \, \, : yA \leq c \}$.
\end{proposition}

\begin{proof}
The polyhedron $P_c$ is the feasible region of $ max \,\,\{ y \cdot b
\,\, : \,\, y A \leq c, \,\, y \in \mathbb R^d\}$, the dual program to
$LP_{A,c}(b)$. The normal fan of $P_c$ is supported on $cone(A)$, i.e.
the union of the normal cones of $P_c$ is $cone(A)$,   
since this is the polar cone of the recession cone $\{y \in
\mathbb R^d : yA \leq 0 \}$ of $P_c$. 
Suppose $b$ is any vector in the interior of a maximal face $cone(A_{\sigma})$
of $\Delta_c$. Then by Proposition~\ref{optbases}, $LP_{A,c}(b)$ has
an optimal solution $x^{\ast}$ with support $\sigma$.
The optimal solution $y$ to the dual of $LP_{A,c}(b)$ satisfies 
$y \cdot a_j = c_j$ for all $j \in \sigma$ and $y \cdot a_j \leq c_j$ 
otherwise, by complementary slackness. Since $\sigma$ is a
maximal face of $\Delta_c$, in fact, $y \cdot a_j < c_j$ for all $j
\not \in \sigma$. This shows that $y$ is unique, and 
$cone(A_{\sigma})$ is contained in the normal cone of $P_c$ at the 
vertex $y$. If $b$ lies in the interior of another maximal face
$cone(A_{\tau})$ then $y'$, the dual optimal solution to $LP_{A,c}(b)$
satisfies $y' \cdot A_{\tau} = c_{\tau}$ and $y' \cdot A_{\bar \tau} <
c_{\bar \tau}$ where $\tau \neq \sigma$. Hence $y'$ is distinct from
$y$ and each maximal cone in $\Delta_c$ lies in a distinct maximal
cone in the normal fan of $P_c$. Since $\Delta_c$ and the normal fan
of $P_c$ have the same support, they must therefore coincide.
\end{proof}

\begin{corollary}
The polyhedron $P_c$ is simple if and only if the regular subdivision
$\Delta_c$ is a triangulation of $cone(A)$. \qed
\end{corollary}

Regular subdivisions were introduced in \cite{GKZ} and have since been
studied from various points of view. They play a central
role in the algebraic study of integer programming (\cite{Stu},
\cite{ST}). We use them here to define {\em group relaxations} of
$IP_{A,c}(b)$.   

A subset $\tau$ of $\{1, \ldots, n\}$ partitions $x = (x_1, \ldots,
x_n)$ as $x_{\tau}$ and $x_{\bar \tau}$ where $x_{\tau}$ consists of
the variables indexed by $\tau$, and $x_{\bar \tau}$ the variables
indexed by the complementary set $\bar \tau$. Similarly, the matrix
$A$ is partitioned as $A = [A_\tau, A_{\bar \tau}]$ and the cost vector
as $c = (c_{\tau}, c_{\bar \tau})$.  If $\sigma$ is a maximal face
of $\Delta_c$ then $A_{\sigma}$ is nonsingular and $Ax = b$ can be
written as $x_{\sigma} = A_{\sigma}^{-1}(b - A_{\bar \sigma}x_{\bar
\sigma})$. Then $c \cdot x = c_{\sigma}(A_{\sigma}^{-1}(b - A_{\bar
\sigma}x_{\bar \sigma})) + c_{\bar \sigma}x_{\bar \sigma} =
c_{\sigma}A_{\sigma}^{-1}b + (c_{\bar \sigma} -
c_{\sigma}A_{\sigma}^{-1}A_{\bar \sigma})x_{\bar \sigma}$. Let
$\tilde{c}_{{\bar \sigma}} := c_{\bar \sigma} -
c_{\sigma}A_{\sigma}^{-1}A_{\bar \sigma}$ and for any face $\tau$ of
$\sigma$, let $\tilde{c}_{{\bar \tau}}$ be the extension of
$\tilde{c}_{{\bar \sigma}}$ to a vector in $\mathbb R^{|\bar \tau|}$ by
adding zeros.

\begin{definition} \label{grouprels}
The group relaxation of the integer program $IP_{A,c}(b)$ with respect  
to the face $\tau$ of $\Delta_c$, denoted as $G^{\tau}(b)$, is the
program $$ min \, \, \{ 
\tilde{c}_{{\bar \tau}} \cdot x_{\bar \tau} \,\, : A_{\tau}x_{\tau} + A_{\bar
  \tau}x_{\bar \tau} = b, \,\, x_{\bar \tau} \geq 0, (x_{\tau}, x_{\bar
  \tau}) \in \mathbb Z^n \}.$$
\end{definition}

The group relaxation $G^{\tau}(b)$ solves $IP_{A,c}(b)$ if its optimal
solution is non-negative. These relaxations contain among them the
usual group relaxations of $IP_{A,c}(b)$ found in the literature. The
program $G^{\sigma}(b)$ where $A_{\sigma}$ is an optimal basis of the
linear relaxation $LP_{A,c}(b)$ is precisely Gomory's group relaxation
of $IP_{A,c}(b)$ \cite{Gom65}.  The set of relaxations $G^{\tau}(b)$
as $\tau$ varies among the subsets of this $\sigma$ are the {\em
extended} group relaxations of $IP_{A,c}(b)$ defined by Wolsey
\cite{Wol}.  Since $\emptyset \in \Delta_c$, $G^{\emptyset}(b) =
IP_{A,c}(b)$ is a group relaxation of $IP_{A,c}(b)$ by
Definition~\ref{grouprels}, and hence $IP_{A,c}(b)$ will certainly be
solved by one of its extended group relaxations. However, it is easy
to construct examples where a group relaxation $G^\tau(b)$ solves
$IP_{A,c}(b)$, but $G^\tau(b)$ is neither Gomory's group relaxation of
$IP_{A,c}(b)$ nor one of its nontrivial extended Wolsey relaxations
(see Example \ref{long-chain}).
Theorem~\ref{bounded} will show that the relaxations in
Definition~\ref{grouprels} are precisely all the bounded group
relaxations of $IP_{A,c}(b)$. Hence this definition considers all the
group relaxations of each integer program in $IP_{A,c}$ that can
possibly solve the program.

For our purposes it is convenient to reformulate $G^{\tau}(b)$ as
follows. Let $B \in \mathbb Z^{n \times (n-d)}$ be any matrix such
that the columns of $B$ generate the $(n-d)$-dimensional lattice
$\mathcal L = \{ x \in \mathbb Z^n : Ax = 0 \} \subset \mathbb Z^n$
and let $u$ be a feasible solution of $IP_{A,c}(b)$. Then \\

$\begin{array}{lll} IP_{A,c}(b) & = & min \,\, \{ c \cdot x \,\,: \,\,
Ax = b, \,\, x \in \mathbb N^n \} \\ 
& = & min \,\, \{ c \cdot x \,\,:
\,\,x \,\equiv \,u \,(mod \,\, \mathcal L), \,\, x \geq 0 \} \\ 
& = & min
\,\, \{ c \cdot x \,\,: \,\,x= u - Bz, \,\, x \geq 0, \,\, z \in
\mathbb Z^{n-d} \} \end{array}$\\  
The last problem is equivalent to  
$min \,\, \{ c \cdot (u - Bz)\,\,: \,\,
Bz \leq u, \,\, z \in \mathbb Z^{n-d} \}$ and hence, $IP_{A,c}(b)$
is equivalent to the problem  
\begin{equation} \label{ip}
min \,\, \{ (-cB)z \,\,: \,\, Bz \leq u, \,\, z \in \mathbb Z^{n-d} \}.
\end{equation}
There is a bijection between the set of feasible solutions of
(\ref{ip}) and the set of feasible solutions of $IP_{A,c}(b)$ via the 
isomorphism $z \mapsto u - Bz$. In particular, $0 \in \mathbb
R^{n-d}$ is feasible for (\ref{ip}) since it is the pre-image of $u$.

Let $\pi_{\tau}$ be the projection map from $\mathbb R^n \rightarrow
\mathbb R^{|\bar \tau|}$ that kills all coordinates indexed by
$\tau$. If $B^{\bar \tau}$ denotes the $|\bar \tau| \times (n-d)$
submatrix of $B$ obtained by deleting the rows indexed by $\tau$, then
we denote by $\mathcal L_{\tau}$ the lattice $\pi_{\tau}(\mathcal L) = 
\{B^{\bar \tau}z \,\, : \,\, z \in \mathbb Z^{n-d} \}$. 
It can be deduced from \cite{SWZ} that the group relaxation
$G^{\tau}(b)$ is equivalent to the {\em lattice program}
$$min \,\, \{ \tilde{c}_{{\bar \tau}} \cdot x_{\bar \tau} \,\,: \,\,
x_{\bar \tau} \equiv \pi_{\tau}(u) \,\, (mod \, \mathcal
L_{\tau}), \,\, x_{\bar \tau} \geq 0 \}$$
which can be reformulated as above to be  
$$
min \,\, \{ (-\tilde{c}_{{\bar \tau}}B^{\bar \tau}) \cdot
z \,\,:\,\, B^{\bar \tau}z \leq \pi_{\tau}(u), \,\, z \in \mathbb
Z^{n-d} \}.
$$

Since $\tilde{c}_{{\bar \tau}} = \pi_{\tau}(c -
c_{\sigma}A_{\sigma}^{-1}A)$ for any maximal face $\sigma$ of
$\Delta_c$ containing $\tau$, and the support of $c -
c_{\sigma}A_{\sigma}^{-1}A$ is contained in $\bar \tau$, we get that 
$\tilde{c}_{{\bar \tau}}B^{\bar \tau} = (c -
c_{\sigma}A_{\sigma}^{-1}A)B = cB$ since $AB = 0$. Hence $G^{\tau}(b)$
is equivalent to 
\begin{equation} \label{gp}
min \,\, \{ (-cB) \cdot z \,\,:\,\, B^{\bar \tau}z \leq
\pi_{\tau}(u), \,\, z \in \mathbb Z^{n-d} \}.
\end{equation}
The feasible solutions to (\ref{ip}) are the lattice points in the
rational polyhedron $P_{u} := \{ z \in \mathbb R^{n-d} : Bz \leq u
\}$ and those to (\ref{gp}) are the lattice points in the relaxation
$P_{u}^{\bar \tau} := \{ z \in \mathbb R^{n-d} : B^{\bar \tau} z
\leq \pi_{\tau}(u) \}$ of $P_{u}$ obtained by deleting the
inequalities indexed by $\tau$. In theory, one could define group
relaxations of $IP_{A,c}(b)$ with respect to any $\tau \subseteq \{1,
\ldots, n\}$. The following result justifies
Definition~\ref{grouprels}.

\begin{theorem} \label{bounded}
The group relaxation $G^{\tau}(b)$ of $IP_{A,c}(b)$ 
has a finite optimal solution if and  
only if $\tau \subseteq \{1, \ldots, n\}$ is a face of $\Delta_c$.
\end{theorem}

\begin{proof}
Since all data are integral it suffices to prove that the linear
relaxation $min \,\, \{ (-cB)z \,\, : \,\, z \in P_{u}^{\bar \tau} \}$
is bounded if and only if $\tau \in \Delta_c$.  

If $\tau$ is a face of $\Delta_c$ then there exists $y \in \mathbb
R^d$ such that $yA_{\tau} = c_{\tau}$ and $yA_{\bar \tau} < c_{\bar
\tau}$. Using the fact that $A_{\tau}B^{\tau} + A_{\bar \tau}B^{\bar
\tau} = 0$ we see that $cB = c_{\tau}B^{\tau} + c_{\bar \tau}B^{\bar
\tau} = yA_{\tau}B^{\tau} + c_{\bar \tau}B^{\bar \tau} = y(-A_{\bar
\tau}B^{\bar \tau}) + c_{\bar \tau}B^{\bar \tau} = (c_{\bar \tau} -
yA_{\bar \tau})B^{\bar \tau}$. This implies that $cB$ is a positive
linear combination of the rows of $B^{\bar \tau}$ since $c_{\bar \tau}
- yA_{\bar \tau} > 0$. Hence $cB$ lies in the polar of $\{z \in
\mathbb R^{n-d} : B^{\bar \tau}z \leq 0 \}$ which is the recession
cone of $P_{u}^{\bar \tau}$ proving that the linear program $min
\,\, \{ (-cB)z \,\, : \,\, z \in P_{u}^{\bar \tau} \}$ is bounded.

The linear program $min \,\, \{ (-cB)z \,\, : \,\, z \in P_{u}^{\bar
\tau} \}$ is feasible since $0$ is a feasible solution. If it is
bounded as well then $min \,\, \{ c_{\tau}x_{\tau} + c_{\bar
\tau}x_{\bar \tau} \,\, : \,\, A_{\tau}x_{\tau} + A_{\bar \tau}x_{\bar
\tau} = b, \, x_{\bar \tau} \geq 0 \}$ is feasible and bounded. Hence
the dual of the latter program $max \,\, \{ y \cdot b \,\, : \,\,
yA_{\tau} = c_{\tau}, \, yA_{\bar \tau} \leq c_{\bar \tau} \}$ is
feasible. This shows that a
superset of $\tau$ is a face of $\Delta_c$ which implies that $\tau
\in \Delta_c$ since $\Delta_c$ is a triangulation.
\end{proof}

The reformulations (\ref{ip}) and (\ref{gp}) imply that $G^{\tau}(b)$
solves $IP_{A,c}(b)$ if and only if both programs have the same
optimal solution $z^{\ast} \in \mathbb Z^{n-d}$. If $G^{\tau}(b)$
solves $IP_{A,c}(b)$ then $G^{\tau'}(b)$ also solves $IP_{A,c}(b)$ for
every $\tau'$ {\em contained in} $\tau$. We say that $\tau \in \Delta_c$ is
{\em associated} to $IP_{A,c}$ if for some $b \in \mathbb N A$,
$G^{\tau}(b)$ solves $IP_{A,c}(b)$ but $G^{\tau'}(b)$ does not for all
faces $\tau' \neq \tau$ of $\Delta_c$ {\em containing} $\tau$.  Several
results about the structure of the subposet of faces of $\Delta_c$
that are associated to $IP_{A,c}$ can be found in \cite{HT2}. For
instance, the associated sets of $IP_{A,c}$ occur in saturated chains
\cite[Theorem 3.1]{HT2}. For a given $b \in \mathbb N A$, the most
easily solved relaxations of $IP_{A,c}(b)$ are those $G^{\sigma}(b)$
where $\sigma$ is a maximal face of $\Delta_c$. We call these
``top-level'' relaxations the {\em Gomory relaxations} of
$IP_{A,c}(b)$.

\begin{definition} \label{gomory-family}
The family of integer programs $IP_{A,c}$ is a 
{\bf Gomory family} if, for {\em every} $b \in \mathbb N A$,
$IP_{A,c}(b)$ is solved by a group relaxation $G^{\sigma}(b)$
where $\sigma$ is a {\em maximal} face of the regular triangulation
$\Delta_c$. 
\end{definition}

Our goal in the rest of this section is to characterize Gomory
families of integer programs (Theorem~\ref{equivconds}). We will
assume from now on that every integer program in $IP_{A,c}$ has a
unique solution which is a stricter notion of genericity of $c$ than
requiring $\Delta_c$ to be a triangulation \cite{ST}. Let $\mathcal
O_c \subseteq \mathbb N^n$ be the set of all optimal solutions of the
programs in $IP_{A,c}$. The set $\mathcal O_c$ is known to be a {\em
  down set} or {\em order ideal} in $\mathbb N^n$, i.e.  $u \in
\mathcal O_c$ and $v \leq u$, $v \in \mathbb N^n$ implies that $v \in
\mathcal O_c$ \cite{Tho}. For a given $A$, there are only finitely
many sets $\mathcal O_c$ as $c$ varies. Two generic cost vectors $c$
and $c'$ are equivalent if $\mathcal O_c = \mathcal O_{c'}$ and all
equivalence classes of generic cost vectors are open full dimensional
cones in $\mathbb R^n$ \cite{ST}. Since $c$ is
generic, $\mathcal O_c$ is in bijection with $\mathbb N A$ via the
linear map $\phi_A : \mathbb N^n \rightarrow \mathbb N A$ where $u
\mapsto Au$.  Let $Q_u := \{ z \in \mathbb R^{n-d} \,\, : \,\, Bz \leq
u, \, (-cB)z \leq 0 \}$ and $Q_u^{\bar \tau} := \{ z \in \mathbb
R^{n-d} \,\, : \,\, B^{\bar \tau} z \leq \pi_{\tau}(u), \, (-cB)z \leq
0 \}$.

\begin{lemma} \label{emptypolys}
(i) A vector $u$ is in  $\mathcal O_c$ if and only
if $Q_{u} \cap \mathbb Z^{n-d} = \{0\}$.\\
(ii) If $u \in \mathcal O_c$, then the group relaxation $G^{\tau}(Au)$
solves the integer program $IP_{A,c}(Au)$ if and only if $Q_{u}^{\bar
\tau} \cap \mathbb Z^{n-d} = \{0\}$.
\end{lemma}

\begin{proof} (i) The lattice point $u$ belongs to $\mathcal O_c$ if
and only if $u$ is the optimal solution to $IP_{A,c}(Au)$ which is
equivalent to $0 \in \mathbb Z^{n-d}$ being the optimal solution to
the reformulation (\ref{ip}) of $IP_{A,c}(Au)$. Since $c$ is generic,
the last statement is equivalent to $Q_{u} \cap \mathbb Z^{n-d} =
\{0\}$. The second statement follows from the fact that (\ref{gp})
solves (\ref{ip}) if and only if they have the same optimal solution.
\end{proof}

By Lemma~\ref{emptypolys}, it is convenient to use the optimal solution
to $IP_{A,c}(b)$ as the vector $u$ in (\ref{ip})
and (\ref{gp}), and we will do so from now on.
For an element $u \in \mathcal O_c$
and a face $\tau$ of $\Delta_c$ let $S(u,\tau)$ be the affine
semigroup $u + \mathbb N (e_i : i \in \tau)$ in $\mathbb
N^n$ where $e_i$ denotes the $i$-th unit vector in $\mathbb R^n$.

\begin{lemma} \label{solvesall}
Suppose $u$ lies in $\mathcal O_c$. If $G^{\tau}(Au)$ solves $IP_{A,c}(Au)$,
then $G^{\tau}(Av)$ solves $IP_{A,c}(Av)$ for all $v \in S(u,\tau)$.
\end{lemma}

\begin{proof}
  If $v \in S(u,\tau)$, then $\pi_{\tau}(u) = \pi_{\tau}(v)$, and this
  implies that $Q_u^{\bar \tau} = Q_v^{\bar \tau}$. If $G^{\tau}(Au)$
  solves $IP_{A,c}(Au)$, then $\{0\} = Q_u^{\bar \tau} \cap \mathbb
  Z^{n-d} = Q_v^{\bar \tau} \cap \mathbb Z^{n-d}$ for all $v \in
  S(u,\tau)$. This implies the result by Lemma~\ref{emptypolys} (ii).
\end{proof}

\begin{proposition} \label{stdpair-grouprel}
The affine semigroup $S(u,\tau)$ is contained
in $\mathcal O_c$ if and only if $G^{\tau}(Au)$ solves $IP_{A,c}(Au)$.
\end{proposition}

\begin{proof} 
Suppose $S(u,\tau) \subseteq \mathcal O_c$. Then for all $v \in
S(u,\tau)$, $Q_v = \{ z \in \mathbb R^{n-d} \,\, : \,\, B^{\tau}z \leq
\pi_{\bar \tau}(v), \, B^{\bar \tau}z \leq \pi_{\tau}(u), \, (-cB)z
\leq 0 \} \cap \mathbb Z^{n-d} = \{0\}$.  Since $\pi_{\bar \tau}(v)$
can be any vector in $\mathbb N^{|\tau|}$, and $Q_u^{\bar \tau}$ is
bounded by Theorem~\ref{bounded}, $Q_u^{\bar \tau} = \{ z \in \mathbb
R^{n-d} \,\, : \,\, B^{\bar \tau}z \leq \pi_{\tau}(u), \, (-cB)z \leq
0 \} \cap \mathbb Z^{n-d} = \{0\}$. Hence, by Lemma~\ref{emptypolys}
(ii), $G^{\tau}(Au)$ solves $IP_{A,c}(Au)$.

Conversely, if $G^{\tau}(Au)$ solves $IP_{A,c}(Au)$, then $\{0\}
= Q_u^{\bar \tau} \cap \mathbb Z^{n-d} = Q_v^{\bar \tau} \cap \mathbb
Z^{n-d}$ for all $v \in S(u,\tau)$.  Since $Q_v^{\bar \tau}$ is a
relaxation of $Q_v$, $Q_v \cap \mathbb Z^{n-d} = \{0\}$ for all $v \in
S(u,\tau)$ and hence by Lemma~\ref{emptypolys} (i), $S(u,\tau)
\subseteq \mathcal O_c$.  
\end{proof}

\begin{definition} 
For $\tau \in \Delta_c$ and $u \in \mathcal O_c$, the pair $(u,\tau)$
is called an {\em admissible pair} of $\mathcal O_c$ if \\ 
\indent (i) $G^{\tau}(Au)$ solves $IP_{A,c}(Au)$ or equivalently,
$S(u,\tau) \subseteq \mathcal O_c$, and \\ 
\indent (ii) the support of $u$ is contained in $\bar \tau$.\\ 
An admissible pair $(u, \tau)$ is a {\em standard pair} of $\mathcal
O_c$ if the affine semigroup $S(u,\tau)$ is not properly contained
in another affine semigroup $S(v,\tau')$ where $(v,\tau')$ is also
an admissible pair of $\mathcal O_c$.
\end{definition}

\begin{definition}\label{stdpolys}
For $u \in \mathbb N^n$ and a face $\tau$ of $\Delta_c$, we say that
the polytope $Q_{u}^{\bar \tau}$ is a {\em standard polytope} of
$IP_{A,c}$ if $Q_{u}^{\bar \tau} \cap \mathbb Z^{n-d} = \{0\}$ and
every relaxation of $Q_{u}^{\bar \tau}$ obtained by removing an
inequality in $B^{\bar \tau}z \leq \pi_{\tau}(u)$ contains a non-zero
lattice point.
\end{definition}

\begin{theorem} \label{stdpairs-equivs}
The following are equivalent:\\
(i) An admissible pair $(u,\tau)$ is a standard pair of $\mathcal
O_c$.\\
(ii) The polytope $Q_{u}^{\bar \tau}$ is a standard polytope of
$IP_{A,c}$.\\ 
(iii) The face $\tau$ of $\Delta_c$ is associated to $IP_{A,c}$.
\end{theorem}

\begin{proof} The proof of $(i) \Leftrightarrow (ii)$ is the content 
of Theorem 2.5 in \cite{HT2}. The equivalence $(i) \Leftrightarrow
(iii)$ follows from the definition of a standard pair and
Lemma~\ref{solvesall}.
\end{proof}

Under the linear map $\phi_A : \mathbb N^n \rightarrow \mathbb N A$
such that $u \mapsto Au$, the affine semigroup $S(u,\tau)$ where
$(u,\tau)$ is a standard pair of $\mathcal O_c$ maps to the affine
semigroup $Au + \mathbb N A_{\tau}$ in $\mathbb N A$. Since every
integer program in $IP_{A,c}$ is solved by one of its group
relaxations, $\mathcal O_c$ is covered by its standard pairs.  We call
this cover and its image in $\mathbb N A$ under $\phi_A$, the {\em
  standard pair decompositions} of $\mathcal O_c$ and $\mathbb N A$
respectively.  Since standard pairs of $\mathcal O_c$ are determined
by the standard polytopes of $IP_{A,c}$, the standard pair
decomposition of $\mathcal O_c$ is unique. The terminology used above
has its origins in \cite{STV} which introduced the {\em standard pair
  decomposition} of a {\em monomial ideal}. The specialization to
integer programming appear in \cite[\S 12.D]{Stu}, \cite{HT1} and
\cite{HT2}. For each $\tau \in \Delta_c$, there are only finitely many
standard pairs of $\mathcal O_c$ that are indexed by $\tau$. Borrowing
terminology from \cite{STV}, we call the number of standard pairs of
the form $(\cdot, \tau)$ the {\em multiplicity} of $\tau$ in $\mathcal
O_c$. The total number of standard pairs is called the {\em arithmetic
  degree} of $\mathcal O_c$. By Theorem~\ref{stdpolys}, multiplicity
of $\tau$ in $\mathcal O_c$ is the number of distinct standard
polytopes of $IP_{A,c}$ indexed by $\bar \tau$ and the arithmetic
degree of $\mathcal O_c$ is the total number of standard polytopes of
$IP_{A,c}$. Lemma~\ref{maxsimplex} shows that the maximal faces of
$\Delta_c$ play a special role in the standard pair decomposition of
$\mathcal O_c$. Part (ii) can also be deduced from the work of Gomory
\cite{Gom65}. 

\begin{lemma} \cite[\S 12.D]{Stu}\label{maxsimplex}
(i) $(0,\sigma)$ is a standard pair of $\mathcal O_c$ if and only if
$\sigma$ is a maximal face of $\Delta_c$.\\ 
(ii) If $\sigma$ is a maximal face of $\Delta_c$, the multiplicity of
$\sigma$ in $\mathcal O_c$ is the index of the sublattice
$\mathbb Z A_{\sigma}$ in $\mathbb Z A$. 
\end{lemma}

The index of $\mathbb Z A_{\sigma}$ in $\mathbb Z A$ is also the
determinant of $A_{\sigma}$ divided by the g.c.d of the
maximal minors of $A$. In contrast to Lemma~\ref{maxsimplex}, if
$\tau$ is a lower dimensional face of $\Delta_c$, it may not index any
standard pair of $\mathcal O_c$. This makes the structure of
$\mathcal O_c$ complicated. (See \cite{HT2} for further results.)
We give an example below.

\begin{example} \label{long-chain}
Consider the following $A \in \mathbb Z^{3 \times 6}$ of rank three:
$$ \left( \begin{array}{cccccc} 
5 & 0 & 0 & 2 & 1 & 0 \\
0 & 5 & 0 & 1 & 4 & 2 \\
0 & 0 & 5 & 2 & 0 & 3
\end{array} \right ).$$
The first three columns of $A$ generate $cone(A)$ which is simplicial. 
If $c = (21, \, 6, \, 1, \, 0, \, 0, \, 0)$ then the regular
triangulation $\Delta_c$ is:
$$\{\{1, 3, 4\}, \, \{1,4,5\}, \,\{2,5,6\}, \, \{3,4,6\}, \,\{4,5,6\}\}.$$ 
The set $\mathcal O_c$ has arithmetic degree 70 which means that
$\mathcal O_c$ has 70 standard pairs which are listed below.
Not all lower dimensional faces of $\Delta_c$ index standard pairs in
this example. 

$$\begin{array}{|c|c|}
\hline
\tau & \text{ standard pairs } (\cdot, \tau) \\
\hline 
\{1, 3, 4\} &  (0, \cdot), \, (e_5, \cdot), \, (e_6, \cdot), \,
  (e_5+e_6, \cdot), \, (2e_6, \cdot) \\  
\hline
\{1,4,5\} & (0, \cdot),\, (e_2, \cdot),\, (e_3, \cdot),\, (e_6, \cdot), \,
 (e_2+e_3, \cdot),\, (2e_2, \cdot),\\ 
     & (3e_2,\cdot),\, (2e_2+e_3, \cdot)\\
\hline
\{2,5,6\}&  (0, \cdot), \, (e_3, \cdot), \, (2e_3, \cdot) \\
\hline
\{3,4,6\}& (0, \cdot), \, (e_5, \cdot), \, (2e_5, \cdot), \, 
           (3e_5, \cdot) \\
\hline
\{4,5,6\}& (0, \cdot), \, (e_3, \cdot), \, (2e_3, \cdot), \, 
           (3e_3, \cdot), \, (4e_3, \cdot) \\
\hline
\{1,4\} & (e_3+2e_5+e_6, \cdot),\,  (2e_3+2e_5+e_6, \cdot), \,
 (2e_3+2e_5, \cdot), \\
        & (2e_3+3e_5, \cdot), \,  (2e_3+4e_5, \cdot) \\
\hline
\{1,5\} & (e_2+e_6, \cdot), \, (2e_2+e_6, \cdot), \, (3e_2+e_6, \cdot) \\
\hline
\{2,5\} & (e_3+e_4, \cdot), \, (e_4, \cdot), \, (2e_4, \cdot) \\
\hline
\{3,4\} & (e_2, \cdot), \,  (e_1+e_2, \cdot), \, (e_1+2e_5, \cdot), \,
          (e_1+2e_5+e_6, \cdot), \, (e_2+e_5, \cdot), \,  \\
\hline
\{3,6\} & (e_2, \cdot), \, (e_2 + e_5, \cdot) \\
\hline
\{4,5\} & (e_2 + 2e_3, \cdot), \, (e_2+3e_3, \cdot), \, 
          (2e_2+2e_3, \cdot), \, (3e_2+e_3, \cdot), \, (4e_2, \cdot) \\
\hline
\{5,6\} & (e_2+3e_3, \cdot) \\
\hline
\{1\} & (e_2+e_3+e_6, \cdot), \, (e_2+e_3+e_5+e_6, \cdot), \,
        (e_2+2e_6, \cdot), \\ 
      & (e_2+e_3+2e_6, \cdot), 
       (2e_2+2e_6, \cdot), \, (e_2+e_3+2e_5+e_6, \cdot)\\
\hline
\{3\} & (e_1+e_2+e_6, \cdot), \, (e_1+e_2+2e_6, \cdot) \\
\hline
\{4\} & (e_1+e_2+2e_3+e_5, \cdot), \,
        (e_1+e_2+2e_3+2e_5, \cdot), \\ 
      & (e_1+e_2+2e_3+3e_5, \cdot), \,
        (e_1+e_2+2e_3+4e_5, \cdot), \\ 
      & (e_1+3e_3+3e_5, \cdot), \,
        (e_1+3e_3+4e_5, \cdot) \\
\hline
\{\emptyset\} & (e_1+e_2+2e_3+e_5+e_6, \cdot), \, 
                (e_1+e_2+2e_3+2e_5+e_6, \cdot), \\  
              & (e_1+2e_2+e_3+e_6, \cdot), \, 
                (e_1+2e_2+e_3+e_5+e_6, \cdot), \\
              & (e_1+2e_2+e_3+2e_5+e_6, \cdot), \, 
                (e_1+2e_2+e_3+2e_6, \cdot), \\
              &  (e_1+3e_2+2e_6, \cdot) \\
\hline 
\end{array}$$

Observe that the integer program $IP_{A,c}(b)$
where $b = A(e_1+e_2+e_3)$
is solved by $G^{\tau}(b)$ with 
$\tau = \{1, 4, 5\}$. 
By Proposition~\ref{optbases}, Gomory's relaxation of $IP_{A,c}(b)$ is
indexed by $\sigma = \{4,5,6\}$ since $b$ lies in the interior of the face
$cone(A_{\sigma})$ of $\Delta_c$. However, neither this relaxation nor
any nontrivial extended relaxation solves $IP_{A,c}(b)$
since the optimal solution $e_1+e_2+e_3$ is not covered by any standard pair
$(\cdot, \tau)$ where $\tau$ is a non-empty subset of $\{4,5,6\}$.
\end{example}

The results stated thus far give characterizations of Gomory families.

\begin{theorem} \label{equivconds}
The following conditions are equivalent:\\
(i) $IP_{A,c}$ is a Gomory family. \\
(ii) The associated sets of $IP_{A,c}$ are precisely the maximal 
faces of $\Delta_c$.\\
(iii) $(\cdot, \tau)$ is a standard pair of $\mathcal O_c$ if and only
if $\tau$ is a maximal face of the regular triangulation
$\Delta_c$.\\ 
(iv) All standard polytopes of $IP_{A,c}$ are simplices.
\end{theorem}

\begin{proof}
The proof follows from Definition~\ref{gomory-family},
Proposition~\ref{stdpair-grouprel} and Theorem~\ref{stdpairs-equivs}.
\end{proof}

If there is a generic cost vector $c$ such that for a triangulation
$\Delta$ of $cone(A)$, $\Delta = \Delta_c$, then we say that $\Delta$
{\em supports} the order ideal $\mathcal O_c$ and the family of
integer programs $IP_{A,c}$. No regular triangulation of the matrix
$A$ in Example~\ref{long-chain} supports a Gomory family. Here is a
matrix with a Gomory family.

\begin{example} \label{example-with-gfamily}
Consider the $3 \times 6$ matrix 
$$A = \left( \begin{array}{cccccc} 
1 & 0 & 1 & 1 & 1 & 1 \\
0 & 1 & 1 & 1 & 2 & 2 \\
0 & 0 & 1 & 2 & 3 & 4
\end{array} \right ).$$
In this case, $cone(A)$ has 14 distinct regular triangulations and 48
distinct sets $\mathcal O_c$ as $c$ varies among all generic cost
vectors. Ten of these triangulations support Gomory families; one for
each triangulation. For
instance, if $c = 
(0,\,0,\,1,\,1,\,0,\,3)$, then
$$\Delta_c = 
\{ \sigma_1 = \{1, 2, 5\}, \,\, \sigma_2 = \{1, 4, 5\}, \,\, 
\sigma_3 = \{2,5,6\}, \,\, \sigma_4 = \{4,5,6\} \}$$
and $IP_{A,c}$ is a Gomory family since the standard pairs of
$\mathcal O_c$ are: $(0, \sigma_1), \, (e_3, \sigma_1), \,  (e_4,
\sigma_1), \,  (0, \sigma_2), \, (0, \sigma_3), \, \text{and} \, (0,
\sigma_4).$ 
\end{example}

\section{Total dual integrality and Gomory families}

We now relate the notion of {\em total dual
integrality} \cite[\S 22]{Sch} to Gomory families. Recall 
that $\mathbb Z A = \mathbb Z^d$ by assumption.

\begin{definition} \label{tdi}
The system $yA \leq c$ is totally dual
integral (TDI) if $LP_{A,c}(b)$ has an integral optimal  
solution for each $b \in cone(A) \cap \mathbb Z^d$. 
\end{definition}

\begin{definition}
The regular triangulation $\Delta_c$ is {\em unimodular} 
if $\mathbb{Z}A_\sigma = \mathbb Z^d$  for every maximal face 
$\sigma \in \Delta_c$. 
\end{definition}

\begin{theorem} \label{TDI-ness} 
The system $yA \leq c$ is TDI if and only if 
the regular triangulation $\Delta_c$ is unimodular.
\end{theorem}

\begin{proof} 
The regular triangulation $\Delta_c$ is the normal fan of $P_c$ by
Proposition~\ref{normalfan}, and it is unimodular if and only if 
$\mathbb Z A_{\sigma} = \mathbb Z^d$ for every
maximal face $\sigma \in \Delta_c$. This is equivalent to saying that
every $b \in cone(A_{\sigma}) \cap \mathbb Z^d$ lies in $\mathbb N
A_{\sigma}$ for every maximal face $\sigma$ of $\Delta_c$. By
Lemma~\ref{optbases}, this happens if and only if $LP_{A,c}(b)$
has an integral optimum for all $b \in cone(A) \cap \mathbb
Z^d$. 
\end{proof}

\begin{corollary} \label{TDI-Gomory}
If $yA \leq c$ is TDI then $IP_{A,c}$ is a Gomory
family.
\end{corollary}

\begin{proof}
  By Lemma~\ref{maxsimplex}, $(0, \sigma)$ is a standard pair of
  $\mathcal O_c$ for every maximal face $\sigma$ of $\Delta_c$.
  Theorem~\ref{TDI-ness} implies that $cone(A_{\sigma})$ is unimodular
  (i.e., $\mathbb Z A_{\sigma} = \mathbb Z^d$), and therefore $\mathbb
  N A_{\sigma} = cone(A_{\sigma}) \cap \mathbb Z^d$ for every maximal
  face $\sigma$ of $\Delta_c$.  Hence the semigroups $\mathbb N
  A_{\sigma}$ arising from the standard pairs $(0,\sigma)$ as $\sigma$
  varies over the maximal faces of $\Delta_c$ cover $\mathbb N A$.
  Therefore the only standard pairs of $\mathcal O_c$ are $(0,\sigma)$
  as $\sigma$ varies over the maximal faces of $\Delta_c$. The result
  then follows from Theorem~\ref{equivconds} (ii).
\end{proof}

When $yA \leq c$ is TDI, the multiplicity of any maximal face $\sigma$
of $\Delta_c$ in $\mathcal O_c$ is one, and all other faces have
multiplicity zero. While this is sufficient for $IP_{A,c}(b)$ to be a
Gomory family, it is far from necessary. TDI-ness guarantees that
$LP_{A,c}(b)$ has an integral optimum for every integral $b$ in
$cone(A)$. In contrast, if $IP_{A,c}(b)$ is a Gomory family, the
linear optima of the programs in $LP_{A,c}$ may not be integral.

If $A$ is {\em unimodular} (i.e., $\mathbb Z A_{\sigma} = \mathbb Z^d$
for every nonsingular maximal submatrix $A_{\sigma}$ of $A$), then the
feasible regions of the linear programs in $LP_{A,c}$ have integral
vertices for each integral $b \in cone(A) \cap \mathbb Z^d$, and $yA
\leq c$ is TDI {\em for all} $c$. Hence if $A$ is unimodular, then
$IP_{A,c}$ is a Gomory family for all generic cost vectors
$c$. However, just as integrality of the optimal solutions of programs
in $LP_{A,c}$ is not necessary for $IP_{A,c}$ to be a Gomory family,
unimodularity of $A$ is not necessary for $IP_{A,c}$ to be a Gomory
family {\em for all} $c$. 

\begin{example}\label{allinishaveembeds}
Consider the seven by twelve integer matrix 
$$ A = \left ( \begin{array}{cccccccccccc}
1&0&0&0&0&0&1&1&1&1&1&0 \\
0&1&0&0&0&0&1&1&0&0&0&1 \\
0&0&1&0&0&0&1&0&1&0&0&1 \\
0&0&0&1&0&0&0&1&0&1&0&0 \\
0&0&0&0&1&0&0&0&1&0&1&0 \\
0&0&0&0&0&1&0&0&0&1&1&1 \\
0&0&0&0&0&0&1&1&1&1&1&1 \\
\end{array} \right )$$ 
of rank seven. The maximal minors of $A$ have absolute values zero,
one and two and hence $A$ is not unimodular. This matrix has $376$
distinct regular triangulations supporting $418$ distinct order ideals
$\mathcal O_c$. In each case, the standard pairs of $\mathcal O_c$ are
indexed by just the maximal simplices of the regular triangulation
$\Delta_c$ that supports it. Hence $IP_{A,c}$ is a Gomory family for
all generic $c$.
\end{example}

\section{$\Delta$-normal matrices}

In Section~3 we saw that unimodularity of $A$ or more generally,
unimodularity of a regular triangulation of $cone(A)$, gives rise to
Gomory families of integer programs. In this section, we identify a
larger set of matrices and cost vectors that give rise to Gomory
families. A common property of unimodular matrices and matrices with a
unimodular triangulation is that they form a {\em Hilbert basis} for
$cone(A)$. In other words, $\mathbb{N}A$ equals $cone(A) \cap
\mathbb{Z}^d$ for such matrices. Borrowing a term from commutative
algebra we make the following definition.

\begin{definition} \label{normal}
A $d \times n$ integer matrix $A$ is {\em normal} if the semigroup
$\mathbb{N}A$ equals $cone(A) \cap \mathbb{Z}^d$.
\end{definition}

We first note that if $A$ is not normal, then $IP_{A,c}$ need not be a
Gomory family for any cost vector $c$. 

\begin{example} \label{nonnormal}
The non-normal matrix $A = \left(\begin{array}{cccc}
1&1&1&1\\0&1&3&4 \end{array}\right)$ gives rise to 10 distinct order
ideals $\mathcal O_c$ supported on its four regular triangulations
$\{\{1,4\}\}, \{\{1,2\},\{2,4\}\},\{\{1,3\},\{3,4\}\}$ and
$\{\{1,2\},\{2,3\},\{3,4\}\}$. Each $\mathcal O_c$ has at least one
standard pair that is indexed by a lower dimensional face of
$\Delta_c$. 
\end{example} 

The matrix in Example~\ref{long-chain} is also not normal and has no
Gomory families. These examples show that normality of $A$ is
necessary for the existence of Gomory families. However, we do not
know at this time whether {\em every} normal matrix $A$ has some
generic cost vector $c$ such that $IP_{A,c}$ is a Gomory family. Our
goal is to show that under certain additional conditions, normal
matrices do give rise to Gomory families.

\begin{definition} A $d \times n$ integer matrix $A$ is
$\Delta$-normal if it has a triangulation $\Delta$  
such that for every maximal face $\sigma \in \Delta$, the columns of $A$ in
$cone(A_{\sigma})$ form a Hilbert basis for $cone(A_{\sigma})$.
\end{definition}

\begin{remark}
  If $A$ is $\Delta$-normal for some triangulation $\Delta$, then it
  is normal. To see this note that every lattice point in $cone(A)$
  lies in $cone(A_{\sigma})$ for some maximal face $\sigma \in
  \Delta$. Since $A$ is $\Delta$-normal, this lattice point also lies
  in the semigroup generated by the columns of $A$ in $cone(A_\sigma)$
  and hence in $\mathbb N A$.
 
  Observe that $A$ is $\Delta$-normal with respect to all the
  unimodular triangulations of $cone(A)$. Hence triangulations
  $\Delta$ with respect to which $A$ is $\Delta$-normal generalize
  unimodular triangulations of $cone(A)$.
\end{remark}

Examples \ref{first} and \ref{second} show that the set of 
matrices where $cone(A)$ has a unimodular triangulation 
is a proper subset of the set of $\Delta$-normal matrices which in
turn is a proper subset of the set of normal matrices.

\begin{example} \label{first}
  Examples of normal matrices with no unimodular triangulations can be
  found in \cite{BoGo} and \cite{FZ}. If $cone(A)$ is simplicial for
  such a matrix, $A$ will be $\Delta$-normal with respect to its
  coarsest (regular) triangulation $\Delta$ consisting of the single
  maximal face with support $cone(A)$.  For
  instance, consider the following example taken from \cite{FZ}:
$$A = \left( \begin{array}{cccccccc} 1 & 0 & 0 & 1 & 1 & 1 & 1 & 1
\\0 & 1 & 0 & 1 & 1 & 2 & 2 & 2 \\0 & 0 & 1 & 1 & 2 & 2 & 3
& 3 \\0 & 0 & 0 & 1 & 2 & 3 & 4 & 5 
\end{array} \right ).$$
Here $cone(A)$ has $77$ regular triangulations and no
unimodular triangulations. Since $cone(A)$ is simplicial, $A$ is
$\Delta$-normal with respect to its coarsest regular triangulation
$\{\{1,2,3,8\}\}$.
\end{example}

\begin{example} \label{second}
There are normal matrices $A$ that are not $\Delta$-normal with respect to
any triangulation of $cone(A)$. To see such an example, consider
the following modification of the matrix in Example~\ref{first} 
that appears in \cite[Example 13.17]{Stu} :
$$A = \left( \begin{array}{ccccccccc} 0 & 1 & 0 & 0 & 1 & 1 & 1 & 1 & 1
\\ 0 & 0 & 1 & 0 & 1 & 1 & 2 & 2 & 2 \\ 0 & 0 & 0 & 1 & 1 & 2 & 2 & 3
& 3 \\ 0 & 0 & 0 & 0 & 1 & 2 & 3 & 4 & 5 \\ 1 & 1 & 1 & 1 & 1 & 1 & 1
& 1 & 1
\end{array} \right ).$$
This matrix is again normal and each of its nine columns generate an 
extreme ray of $cone(A)$. Hence the only way for this matrix to be
$\Delta$-normal for some $\Delta$ would be if $\Delta$ is a 
unimodular triangulation of $cone(A)$. However, this $cone(A)$ has no
unimodular triangulations.
\end{example}

\begin{theorem} \label{specialinitial}
If $A$ is $\Delta$-normal for some regular triangulation $\Delta$ then
there exists a generic cost vector $c \in \mathbb Z^n$ such that
$\Delta = \Delta_c$ and $IP_{A,c}$ is a Gomory family.
\end{theorem}

\begin{proof}
Without loss of generality we can assume that the columns of $A$ in
$cone(A_{\sigma})$ form a minimal Hilbert basis of this cone for 
any maximal face $\sigma$ of $\Delta$. If there were a
redundant element, the smaller matrix obtained by removing this
column from $A$ would still be $\Delta$-normal.

For a maximal face $\sigma \in \Delta$, let $\sigma_{in} \subset
\{1,\ldots,n\}$ be the set of indices of all columns of $A$ lying in
$cone(A_{\sigma})$ that are different from the columns of
$A_{\sigma}$. Suppose $a_{i_1}, \ldots, a_{i_k}$ are the columns of
$A$ that generate the one dimensional faces of $\Delta$, and $c' \in
\mathbb R^n$ a cost vector such that $\Delta = \Delta_{c'}$. We modify
$c'$ to obtain a new cost vector $c \in \mathbb R^n$ such that $\Delta
= \Delta_c$ as follows. For $j = 1, \ldots, k$, let $c_{i_j} :=
c'_{i_j}$. If $j \in \sigma_{in}$ for some maximal face $\sigma \in
\Delta$, then $a_j =\sum_{i \in \sigma} \lambda_ia_i$, $0 \leq
\lambda_i < 1$ and we define $c_j := \sum_{i \in \sigma}
\lambda_ic_i$. Hence, for all $j \in \sigma_{in}$, $(a_j, c_j) \in
\mathbb R^{d+1}$ lies in $C_{\sigma} := cone((a_i,c_i) : i \in \sigma)
= cone((a_i,{c'}_i) : i \in \sigma)$ which was a facet of $C =
cone((a_i,{c'}_i) : i = 1, \ldots, n)$. If $y \in \mathbb R^d$ is a
vector as in Definition~\ref{regtriang} showing that $\sigma$ is a
maximal face of $\Delta_{c'}$ then $y \cdot a_i = c_i$ for all $i \in
\sigma \cup \sigma_{in}$ and $y \cdot a_j < c_j$ otherwise. Since
$cone(A_{\sigma}) = cone(A_{\sigma \cup \sigma_{in}})$, we conclude
that $cone(A_{\sigma})$ is a maximal face of $\Delta_c$.

If $b \in \mathbb{N}A$ lies in $cone(A_\sigma)$ for a maximal face
$\sigma \in \Delta_c$, then $IP_{A,c}(b)$ has at least one feasible
solution $u$ with support in $\sigma \cup \sigma_{in}$ since $A$ is
$\Delta$-normal. Further, $(b, c \cdot u) = (Au,c \cdot u)$ lies in
$C_{\sigma}$ and all feasible solutions of $IP_{A,c}(b)$ with support
in $\sigma \cup \sigma_{in}$ have the same cost value by
construction. Suppose $v \in \mathbb{N}^n$ is any feasible solution of
$IP_{A,c}(b)$ with support not in $\sigma \cup \sigma_{in}$. Then $c
\cdot u < c \cdot v$ since $(a_i,c_i) \in C_{\sigma}$ if and only if
$i \in \sigma \cup \sigma_{in}$ and $C_{\sigma}$ is a facet in the
{\em lower envelope} of $C$. Hence the optimal solutions of
$IP_{A,c}(b)$ are precisely those feasible solutions with support in
$\sigma \cup \sigma_{in}$. The vector $b$ can be expressed as $b = b'
+ \sum_{i \in \sigma} z_ia_i$ where $z_i \in \mathbb N $ are unique
and $b' \in \{\sum_{i \in \sigma} \lambda_i a_i \, : \, 0 \leq
\lambda_i < 1\} \cap \mathbb{Z}^d$ is also unique. The vector $b' =
\sum_{j \in \sigma_{in}} r_ja_j$ where $r_j \in \mathbb N$. Setting
$u_i = z_i$ for all $i \in \sigma$, $u_j = r_j$ for all $j \in
\sigma_{in}$ and $u_k = 0$ otherwise, we obtain all feasible solutions
$u$ of $IP_{A,c}(b)$ with support in $\sigma \cup \sigma_{in}$.

If there is more than one such feasible solution, then $c$ is not
generic. In this case, we can perturb $c$ to a generic cost vector
$c'' = c + \epsilon \omega$ by choosing $1 \gg \epsilon > 0$, $\omega_j
\ll 0$ whenever $j = i_1, \ldots, i_k$ and $\omega_j = 0$ otherwise.
Suppose $u_1, \ldots, u_t$ are the optimal solutions of the integer
programs $IP_{A,c''}(b')$ where $b' \in \{\sum_{i \in \sigma}
\lambda_i a_i \, : \, 0 \leq \lambda_i < 1\} \cap \mathbb{Z}^d$. (Note
that $t = |\{\sum_{i \in \sigma} \lambda_i a_i \, : \, 0 \leq
\lambda_i < 1\} \cap \mathbb{Z}^d|$ is the index of $\mathbb Z
A_{\sigma}$ in $\mathbb Z A$.)  The support of each such $u_i$ is
contained in $\sigma_{in}$.  For any $b \in cone(A_{\sigma}) \cap
\mathbb Z^d$, the optimal solution of $IP_{A,c''}(b)$ is hence $u =
u_i + z$ for some $i \in \{1, \ldots, t\}$ and $z \in \mathbb
N^n$ with support in $\sigma$. This shows that $\mathbb N A$
is covered by the affine semigroups $\phi_A(S(u_i,\sigma))$ where
$\sigma$ is a maximal face of $\Delta$ and $u_i$ as above for each
$\sigma$. By construction, the corresponding admissible pairs
$(u_i,\sigma)$ are all standard for $\mathcal O_{c''}$. Since all data 
is integral, $c'' \in \mathbb Q^n$ and hence can be scaled to lie 
in $\mathbb Z^n$. Renaming $c''$ as $c$, we conclude that
$IP_{A,c}$ is a Gomory family. 
\end{proof}

\begin{corollary} Let $A$ be any normal matrix such that $cone(A)$ is
  simplicial, and let $\Delta$ be the coarsest triangulation 
  whose single maximal face has support $cone(A)$. Then there exists a
  cost vector $c \in \mathbb Z^n$ such that $\Delta = \Delta_c$ and 
  $IP_{A,c}$ is a Gomory family.  
\end{corollary}

\begin{example} \label{non-delta-Gomory}
Consider the normal matrix in Example \ref{example-with-gfamily}.
Here $cone(A)$ is generated by the first, second and sixth columns of
$A$ and hence $A$ is $\Delta$-normal with respect to the regular
triangulation $\{\{1, \, 2, \, 6\}\}$. There are 13 distinct sets
$\mathcal O_c$ supported on $\Delta$. Among the 13 corresponding
families of integer programs, only one is a Gomory family. A
representative cost vector for this $IP_{A,c}$ is
$c=(0,0,4,4,1,0)$. The standard pair decomposition of $\mathcal O_c$
is the one constructed in Theorem \ref{specialinitial}. The affine
semigroups $S(\cdot, \sigma)$ from this decomposition are:
$$ S(0, \sigma), \,\, S(e_3, \sigma), \,\, S(e_4, \sigma), \,\,
\text{and} \,\, S(e_5, \sigma).$$ 
Note that $A$ is not $\Delta$-normal
with respect to the regular triangulation supporting the Gomory family
$IP_{A,c}$ in Example \ref{example-with-gfamily}. The columns of $A$
in $cone(A_{\sigma_1})$ are the columns of $A_{\sigma_1}$ and $A_3$.
The vector $(1, \, 2, \, 2)$ is in the minimal Hilbert basis of 
$cone(A_{\sigma_1})$ but is not a column of $A$. This example shows
that a regular triangulation $\Delta$ of $cone(A)$ can support a
Gomory family even if $A$ is not $\Delta$-normal. The Gomory families 
in Theorem~\ref{specialinitial} have a very special standard pair 
decomposition.
\end{example}


\section{Hilbert covers and Gomory families }

The results in the previous section lead  to the following problem.

\begin{problem} \label{Gomory-question}
If $A \in \mathbb Z^{d \times n}$ is a normal matrix, does there
exist a generic cost vector $c \in \mathbb Z^n$ such that $IP_{A,c}$
is a Gomory family? 
\end{problem}

We do not know the answer to this question. However, in this section
we answer a stronger version of this question for small values of $d$
and state our observations for general $d$. We begin with the
following result. 

\begin{theorem} \label{unimodular-for-three}
If $A \in \mathbb Z^{d \times n}$ is a normal matrix and $d \leq 3$,
then there exists a generic cost vector $c \in \mathbb Z^n$ such that
$IP_{A,c}$ is a Gomory family.
\end{theorem}

\begin{proof}
It is known that if $d \leq 3$ then $cone(A)$ has a regular unimodular 
triangulation $\Delta_c$ \cite{Seb}. The result then follows from 
Corollary \ref{TDI-Gomory}.
\end{proof}

Before we proceed, we rephrase Problem~\ref{Gomory-question} in terms
of covering properties of $cone(A)$ and $\mathbb N A$ along the lines
of \cite{BoGo}, \cite{BG}, \cite{BGHMW}, \cite{FZ} and \cite{Seb}. To
obtain the same set up as in these papers we assume in this section
that $A$ is normal and the columns of $A$ form the unique minimal
Hilbert basis of $cone(A)$. Using the terminology in \cite{BG}, the
{\em free Hilbert cover} problem asks whether there exists a covering
of $\mathbb N A$ by semigroups $\mathbb{N}A_\tau$ where the columns of
$A_\tau$ are linearly independent. The {\em unimodular Hilbert cover}
problem asks whether $cone(A)$ can be covered by full 
dimensional unimodular subcones
$cone(A_{\tau})$ (i.e., $\mathbb Z A_{\tau} = \mathbb Z^d$), while the
stronger {\em unimodular Hilbert partition} problem asks whether
$cone(A)$ has a unimodular triangulation. (Note that if $cone(A)$ has
a unimodular Hilbert cover or partition using subcones
$cone(A_{\tau})$, then $\mathbb N A$ is covered by the semigroups
$\mathbb N A_{\tau}$.)  All these problems have positive answers if $d
\leq 3$ since $cone(A)$ admits a unimodular Hilbert partition in this
case \cite{BoGo}, \cite{Seb}. Normal matrices (with $d = 4$) such that
$cone(A)$ has no unimodular Hilbert partition can be found in
\cite{BoGo} and \cite{FZ}. Examples (with $d = 6$) that admit no free
Hilbert cover and hence no unimodular Hilbert cover can be found in
\cite{BG} and \cite{BGHMW}.

When $yA \leq c$ is TDI, the standard pair decomposition of $\mathbb N A$
induced by $c$ gives a unimodular Hilbert partition of $cone(A)$ by
Theorem~\ref{TDI-ness}. An important difference between
Problem~\ref{Gomory-question} and the Hilbert cover problems is that
{\em affine} semigroups cannot be used in Hilbert covers. Moreover, 
affine semigroups that are allowed in standard pair decompositions
come from integer programming. If there are no restrictions on the
affine semigroups that can be used in a cover, $\mathbb N A$ can
always be covered by full dimensional affine semigroups: for any
triangulation $\Delta$ of $cone(A)$ with maximal subcones 
$cone(A_{\sigma})$, the affine semigroups $b + \mathbb{N}A_\sigma$
cover $\mathbb N A$ as $b$ varies in $\{ \sum_{i \in \sigma} \lambda_i
a_i: \,\, 0 \leq \lambda_i < 1\} \cap \mathbb{Z}^d$ and $\sigma$
varies among the maximal faces of the triangulation. A {\em partition}
of $\mathbb N A$ derived from this idea can be found in 
\cite[Theorem 5.2]{Sta82}.

In order to state our main theorem, we recall the notion of 
{\em supernormality} which was introduced in \cite{HMS}.

\begin{definition} \label{supernormal}
A matrix $A \in \mathbb Z^{d \times n}$ is {\em supernormal}
if for every submatrix $A'$ of $A$, the columns of $A$ that 
lie in $cone(A')$ form a Hilbert basis for $cone(A')$.  
\end{definition}

\begin{proposition} \label{equivalence}
For $A \in \mathbb Z^{d \times n}$, the following are equivalent:
\begin{enumerate}
\item[(i)] $A$ is supernormal,
\item[(ii)] $A$ is $\Delta$-normal for every regular triangulation
$\Delta$ of $cone(A)$,
\item[(iii)] Every triangulation of $cone(A)$ in which all columns of 
$A$ generate one dimensional faces is unimodular.
\end{enumerate}
\end{proposition}

\begin{proof} The equivalence of (i) and (iii) was established in 
\cite[Proposition 3.1]{HMS}. Definition \ref{supernormal} shows that
(i) $\Rightarrow$ (ii). Hence we just need to show that (ii)
$\Rightarrow$ (i). Suppose that $A$ is $\Delta$-normal for every
regular triangulation of $cone(A)$. In order to show that $A$ is supernormal
we only need to check submatrices $A'$ where the dimension of
$cone(A')$ is $d$.  Choose a cost vector $c$ with $c_i \gg 0$ if the
$i$-th column of $A$ does not generate an extreme ray of $cone(A')$,
and $c_i = 0$ otherwise. This gives a polyhedral subdivision of
$cone(A)$ in which $cone(A')$ is a maximal face. There are
standard procedures that will refine this subdivision to a regular
triangulation $\Delta$ of $cone(A)$. Let $T$ be the set of maximal faces
$\sigma$ of $\Delta$ such that $cone(A_{\sigma})$ lies in
$cone(A')$. Since $A$ is $\Delta$-normal, the columns of $A$ that lie
in $cone(A_{\sigma})$ form a Hilbert basis for $cone(A_{\sigma})$ for
each $\sigma \in T$. However, since their union is the set of columns
of $A$ that lie in $cone(A')$, this union forms a Hilbert basis for
$cone(A')$.
\end{proof}

It is easy to catalog all $\Delta$-normal and supernormal matrices, of
the type considered in this paper, for small values of $d$. We say
that the matrix $A$ is {\em graded} if its columns span an affine
hyperplane in $\mathbb R^d$. If $d=1$, $cone(A)$ has $n$
triangulations $\{\{i\}\}$ each of which has the unique maximal
subcone $cone(A_i)$ whose support is $cone(A)$. If 
we assume that $a_1 \leq a_2 \leq \cdots \leq a_n$, then $A$ is normal
if and only if either $a_1 = 1$, or $a_n = -1$. Also, $A$ is normal if
and only if it is supernormal. If $d=2$ and the columns of $A$ are
ordered counterclockwise around the origin, then $A$ is normal if and
only if $det(a_i, a_{i+1}) = 1$ for all $i = 1, \ldots, n-1$. Such an
$A$ is supernormal since it is $\Delta$-normal for every triangulation
$\Delta$ --- the Hilbert basis of a maximal subcone of $\Delta$ is
precisely the set of columns of $A$ in that subcone. If $d = 3$ then
as mentioned before, $cone(A)$ has a unimodular triangulation with
respect to which $A$ is $\Delta$-normal. However, not every such $A$
needs to be supernormal: we saw that the matrix in Example
\ref{example-with-gfamily} is not $\Delta$-normal for the $\Delta$ 
supporting the Gomory family in that example.
If $d=3$ and $A$ is graded, then
without loss of generality we can assume that the columns of $A$ span
the hyperplane $x_1 = 1$. If $A$ is normal as well, then its columns
are precisely all the lattice points in the convex hull of
$A$. Conversely, every graded normal $A$ with $d=3$ arises this way ---
its columns are all the lattice points in a polygon in $\mathbb R^2$
with integer vertices. In particular, every triangulation of $cone(A)$
that uses all the columns of $A$ is unimodular. Hence, by Proposition
\ref{equivalence}, $A$ is supernormal, and therefore $\Delta$-normal
for any triangulation of $A$. 

\begin{theorem}\label{smalld}
Let $A \in \mathbb Z^{d \times n}$ be a normal matrix of rank $d$.

\begin{enumerate}

\item[(i)] If $d = 1,2$ or $A$ is graded and $d = 3$, every
  regular triangulation of $cone(A)$ supports at least one
  Gomory family.

\item[(ii)] If $d = 2$ and $A$ is graded, every regular
  triangulation of $cone(A)$ supports exactly one Gomory family.

\item [(iii)] If $d = 3$ and $A$ is not graded, or if 
  $d = 4$ and $A$ is graded, then not all regular
  triangulations of $cone(A)$ may support a Gomory family. In
  particular, $A$ may not be $\Delta$-normal with respect to every
  regular triangulation. 
\end{enumerate}
\end{theorem}

\begin{proof}
(i) If $d=1,2$ or $A$ is graded and $d=3$, $A$ is supernormal and
hence by Proposition~\ref{equivalence} and
Theorem~\ref{specialinitial}, every regular triangulation of $cone(A)$
supports at least one Gomory family.

\smallskip
\noindent
(ii) If $d=2$ and $A$ is graded, then we may assume that
$$A = \left( \begin{array}{ccccc}
    1 & 1 & 1 & \ldots & 1 \\
    0 & 1 & 2 & \ldots & n-1 \end{array} \right).$$
In this case, $A$
is supernormal and hence every regular triangulation $\Delta$ of
$cone(A)$ supports a Gomory family by Theorem~\ref{specialinitial}.
Suppose the maximal cones of $\Delta$, in counter-clockwise order, are
$C_1, \ldots, C_r$. Assume the columns of $A$ are labeled
such that $C_i = cone(a_{i-1},a_i)$ for $i = 1, \ldots, r$, and 
the columns of $A$ in the interior of $C_i$ are labeled in
counter-clockwise order as $b_{i1}, \ldots, b_{ik_i}$. Hence the $n$
columns of $A$ from left to right are:
$$a_0, b_{11}, \cdots, b_{1k_1}, a_{1}, b_{21}, \cdots, a_{r-1},
b_{r1}, \cdots, b_{rk_r}, a_{r}.$$
Indexing the columns of $A$ by
their labels, the maximal faces of $\Delta$ are $\sigma_i = \{
i-1,i\}$ for $i = 1, \ldots, r$. Let $e_i$ be the unit vector of
$\mathbb R^n$ indexed by the true column index of $a_i$ in $A$ and
$e_{ij}$ be the unit vector of $\mathbb R^n$ indexed by the true
column index of $b_{ij}$ in $A$.  Since the columns of $A$ form a
minimal Hilbert basis of $cone(A)$, $e_i$ is the unique solution to
$IP_{A,c}(a_i)$ for all $c$ and $e_{ij}$ is the unique solution to
$IP_{A,c}(b_{ij})$ for all $c$.  Hence the standard pairs of
Theorem~\ref{specialinitial} are $(0,\sigma_i)$ and
$(e_{ij},\sigma_i)$ for $i=1, \ldots, r$ and $j = 1, \ldots, k_i$.

Suppose $\Delta$ supports a second Gomory family $IP_{A,\omega}$.
Then every standard pair of ${\mathcal O}_w$ is also of the form
$(\cdot, \sigma_i)$ for $\sigma_i \in \Delta$, and $r$ of them are
$(0, \sigma_i)$ for $i = 1, \ldots, r$. The remaining standard pairs
are of the form $(e_{ij}, \sigma_k)$. To see this, consider the
semigroups in $\mathbb N A$ arising from the standard pairs of
$\mathcal O_w$. The total number of standard pairs of $\mathcal
O_c$ and $\mathcal O_w$ are the same. Since the columns of $A$ all lie 
on $x_1 =1$, no two $b_{ij}$s can be covered by a semigroup
coming from the same standard pair and none of them are covered by a
semigroup $(0, \sigma_i)$. We show that if $(e_{ij},
\sigma_k)$ is a standard pair of $\mathcal O_w$ then $k = i$ and thus
$\mathcal O_w = \mathcal O_c$.

If $r=1$, the standard pairs of $\mathcal O_w$ are $(0,\sigma_1),
(e_{11}, \sigma_1), \ldots, (e_{1k_1}, \sigma_1)$ as in
Theorem~\ref{specialinitial}. If $r > 1$, consider the last cone $C_r
= cone(a_{r-1}, a_r)$. If $a_{r-1}$ is the second to last column of $A$,
then $C_r$ is unimodular and the semigroup from $(0,\sigma_r)$ covers
$C_r \cap \mathbb Z^2$. The subcomplex comprised of $C_1, \ldots,
C_{r-1}$ is a regular triangulation $\Delta'$ of $cone(A')$ where $A'$
is obtained by dropping the last column of $A$. Since $A'$ is a normal
graded matrix with $d=2$ and $\Delta'$ has less than $r$ maximal
cones, the standard pairs supported on $\Delta'$ are as in
Theorem~\ref{specialinitial} by induction. If $a_{r-1}$ is not the
second to last column of $A$ then $b_{rk_r}$, the second to last column of
$A$ is in the Hilbert basis of $C_r$ but is not a generator of $C_r$.
So ${\mathcal O}_w$ has a standard pair of the form $(e_{rk_r},
\sigma_i)$.  If $\sigma_i \neq \sigma_r$, then the lattice point
$b_{rk_r} + a_r$ cannot be covered by the semigroup from this or any
other standard pair of $\mathcal O_w$. Hence $\sigma_i = \sigma_r$. By
a similar argument, the remaining standard pairs indexed by $\sigma_r$
are $(e_{r (k_r-1)}, \sigma_r), \ldots,(e_{r1},\sigma_r)$ along with
$(0,\sigma_r)$. These are precisely the standard pairs of $\mathcal
O_c$ indexed by $\sigma_r$. Again we are reduced to considering the
subcomplex comprised of $C_1, \ldots, C_{r-1}$ and by induction, the
remaining standard pairs of $\mathcal O_w$ are as in
Theorem~\ref{specialinitial}.

\smallskip
\noindent
(iii) The 
$3 \times 6$ normal matrix $A$  of Example \ref{example-with-gfamily} 
has 10 distinct Gomory families
supported on 10 out of the 14 regular triangulations of $cone(A)$.
Furthermore, the normal matrix $$A = \left( \begin{array}{ccccccc}
1& 1& 1& 1& 1& 1& 1 \\
1& 0& 1& 1& 1& 1& 0\\
0& 1& 2& 2& 1& 1& 0\\
0& 0& 4& 3& 2& 1& 0 \end{array} \right)$$ has 11 distinct Gomory 
families supported on 11 out of the 19 regular triangulations of
$cone(A)$.
\end{proof}

\section{Computations and Commutative Algebra}

As mentioned in the introduction, the computations
in this paper were done using the connections of this material to
commutative algebra (see \cite{Stu}, \cite{ST} and \cite{Tho}). In
this section we give a brief description of our methods. 

Our general feeling is that Problem~\ref{Gomory-question} has a
negative answer. To check whether a matrix $A$ has a Gomory family, we
first need to compute all the distinct sets of optimal solutions
$\mathcal O_c$ (to the programs in $IP_{A,c}$) that arise as $c$
varies among the generic cost vectors with respect to $A$. As
mentioned in Section~2, there are only finitely many such sets for a
fixed $A$. To check whether $IP_{A,c}$ is a Gomory family, we need to
compute the standard pair decomposition of $\mathcal O_c$: $IP_{A,c}$
is a Gomory family if and only if all the standard pairs of $\mathcal
O_c$ are indexed by maximal faces of $\Delta_c$.

A {\em monomial} $x^u$ in the {\em polynomial ring } $S :=
\mathbb{Q}[x_1, \ldots, x_n]$ is the product 
$x^u = x_1^{u_1}\cdots x_n^{u_n}$ where $u = (u_1, \ldots, u_n) \in
\mathbb N^n$. An ideal $M$ in $S$ is called a {\em monomial ideal} if
it is generated by monomials. Since every ideal in $S$ is finitely
generated, $M = \langle x^{v_1}, \ldots, x^{v_t} \rangle$ for a set of
minimal generators $x^{v_1}, \ldots, x^{v_t}$. The {\em toric ideal}
of $A$, denoted as $I_A$ is the {\em binomial ideal} in $S$ defined
as:
$$ \langle x^u - x^v : \,\, u, v \in \mathbb{N}^n \text{ and } Au = Av
\rangle. $$ The {\em cost} of a monomial $x^u$ with respect to a cost
vector $c \in \mathbb R^n$ is the dot product $c \cdot u$ and the {\em
initial term} of a polynomial $f = \sum \lambda_u x^u \in S$ is
the sum of all terms in $f$ of highest cost. For any ideal $I \subset
S$, the {\em initial ideal} of $I$ with respect to $c$, denoted as
$in_c(I)$, is the ideal generated by all the initial terms of all
polynomials in $I$. These concepts come from the theory of {\em
Gr\"obner bases} for polynomial ideals \cite{CLO}. The toric ideal
$I_A$ provides the algebraic link between integer programming and
Gr\"obner basis theory. For an introduction to this connection see
\cite{Stu}. 

\begin{proposition}
(i) A cost vector $c \in \mathbb R^n$ is generic for $A$ if and only
if the initial ideal $in_c(I_A)$ is a monomial ideal.\\
(ii) For a generic $c$, the vector $u$ belongs to $\mathcal O_c$ if
and only if $x^u$ is not in the intial ideal $in_c(I_A)$. \qed
\end{proposition}

There are only finitely many distinct initial ideals for a 
polynomial ideal \cite{Stu}, and hence there are only
finitely many distinct sets $\mathcal O_c$ as $c$ varies among the 
generic cost vectors. We say that two cost vectors $c$ and $c'$ in
$\mathbb R^n$ are {\em equivalent} if $in_c(I_A) = in_{c'}(I_A)$.

\begin{theorem} \cite[Theorem 3.10]{ST}
(i) Each equivalence class of generic cost vectors form an open full
dimensional polyhedral cone in $\mathbb R^n$. The closure of this cone
is called the {\em Gr\"obner cone} of $c$ where $c$ is any vector 
in the interior of the cone.\\
(ii) The collection of all Gr\"obner cones of $A$ form a complete 
polyhedral fan in $\mathbb R^n$ called the {\em Gr\"obner fan} of $A$.\\
(iii) The Gr\"obner fan of $A$ is the normal fan of an 
$(n-d)$-dimensional polytope called the {\em state polytope} of $A$. \qed
\end{theorem}

By the above results, $\mathcal O_c$ can be computed implicitly by
computing the monomial ideal $in_c(I_A)$. This can be done using a
computer algebra package like Macaulay 2 \cite{GrS}. In order to find
all initial ideals of $I_A$, we use the software package TiGERS
\cite{HuTh} for enumerating the vertices of the state polytope of $A$.
At each vertex, TiGERS returns the initial ideal induced by a vector
in the interior of the normal cone at that vertex.  The standard pair
decomposition of the set of monomials outside a monomial ideal
described in terms of its minimal generators can be calculated using
Macaulay 2. See the chapter {\em Monomial Ideals} in \cite{HoSm}.

To obtain a normal matrix of the type discussed in this paper, it
suffices to start with an arbitrary set of vectors $a_1, \ldots, a_p
\in \mathbb Z^d$ such that $cone(a_1, \ldots, a_p)$ is pointed and
full dimensional and then to compute the Hilbert basis of $cone(a_1,
\ldots, a_p)$. The elements in the Hilbert basis form the columns of a
normal matrix. We used the package Normaliz by Bruns and Koch
\cite{BK} to compute Hilbert bases.

The regular triangulation $\Delta_c$ of $cone(A)$ is a pure
$d$-dimensional complex of cones and hence every maximal face $\sigma
\in \Delta_c$ has cardinality $d$. Hence we get the following algorithm.

\begin{algorithm} \label{check} 
{\em How to check if $A \in \mathbb Z^{d \times n}$
gives a Gomory family.}\\
(i) Compute all the initial ideals of the toric ideal $I_A$ using
TiGERS.\\
(ii) For each initial ideal $in_c(I_A)$:\\
\indent (a) Use Macaulay 2 to find its standard pairs.\\
\indent (b) $IP_{A,c}$ is a Gomory family if and only if 
every set $\tau$ in a standard pair $(\cdot, \tau)$ has
cardinality $d$.
\end{algorithm}

A monomial $x^u$ is {\em square-free} if $u \in \{0,1\}^n$ and a
monomial ideal $M$ is {\em square-free} if all its minimal generators
are square-free. The {\em radical} of a monomial ideal $M$ in $S$ is
the ideal $\sqrt{M} := \langle f : f^r \in M \,\, \text {for some}
\,\,r \in \mathbb N \rangle$. The radical $\sqrt{M}$ is a square-free
monomial ideal. The {\em Stanley Reisner} ideal of the regular
triangulation $\Delta_c$ is the square-free monomial ideal
$$\langle \Pi_{i \in \tau}x_i : \tau \,\, \text {is a minimal non-face of} 
\,\, \Delta_c \rangle.$$

\begin{theorem} \cite[Theorem 8.3]{Stu}
The Stanley-Reisner ideal of the regular triangulation $\Delta_c$ is
the radical of $in_c(I_A)$. \qed
\end{theorem}

If two distinct initial ideals $in_c(I_A)$ and $in_{c'}(I_A)$ have the
same radical, then $\Delta_c = \Delta_{c'}$ and we say that $\Delta_c$
supports these initial ideals. Several initial ideals of $I_A$ may
have the same radical. Using the above theorem, the initial ideals of
$I_A$ output by TiGERS can be grouped according to their radicals or
equivalently, the regular triangulations supporting them. This, in
combination with Algorithm~\ref{check}, allows us to check whether a
regular triangulation supports a Gomory family.

Recall from Theorem~\ref{TDI-ness} that $yA \leq c$ is TDI if
and only if $\Delta_c$ is unimodular. Using the following result,
we obtain an algebraic check for TDI-ness of $yA \leq c$.

\begin{theorem} \cite[Corollary 8.9]{Stu}
The regular triangulation $\Delta_c$ is unimodular if and only if 
the initial ideal $in_c(I_A)$ is square-free. \qed
\end{theorem}

%
%

\bibliography{refs}
\bibliographystyle{plain}

\end{document}